\documentstyle[12pt,leqno]{article}
\title{
\bf Deformation of fillable $CR$ structures}
\author{
  Jih-Hsin Cheng
\thanks{
Research supported in part by National Science
Council grant NSC 90-2115-M-001-001 (R.O.C.) and
a grant of Academia Sinica, Taipei, Taiwan. 
}}
\date{}

\begin{document}
\maketitle

\begin{abstract}
We study the fillability (or embeddability) of $CR$ structures
under the gauge-fixed Cartan flow. We prove that if the initial $CR$
structure
is fillable with nowhere vanishing Tanaka-Webster curvature and free
torsion, then it keeps having the same property after a short time. In the
Appendix, we show the uniqueness of the solution to the gauge-fixed Cartan
flow.
\ \ \ \ \ \ \ \ \ \ \ \ \ \ \ \ \

\end{abstract}

\ \ \ \ \ \ \ \ \ \ \ \ \ \ \ \ \

\bigskip

\paragraph{1. Introduction}

In [CL1], we study an evolution equation for $CR$ structures $J_{(t)}$on a
closed (compact with no boundary) contact 3-manifold $(M,\xi)$ according to
their Cartan (curvature) tensor $Q_{J_{(t)}}$(see also \S2):

\bigskip

$(1.1)$ \ \ \ \ \ \ $\frac{\partial J_{(t)}}{\partial t}=2Q_{J_{(t)}}.$

\bigskip

We will often call this evolution equation (1.1) the Cartan flow.
Since the equation $(1.1)$ is invariant under a big symmetry group, namely,
the contact diffeomorphisms, we add a gauge-fixing term on the right-hand
side
to break the symmetry. The gauge-fixed (called ''regularized'' in [CL1])
Cartan flow reads as follows:

\bigskip

$(1.2)$ \ \ \ \ $\frac{\partial J_{(t)}}{\partial t}=2Q_{J_{(t)}}-\frac{1}%
{6}D_{J_{(t)}}F_{J_{(t)}}K$

\bigskip

\noindent (see [CL1] or \S\ 2 for the meaning of notations). In this paper,
we
investigate the fillability of $CR$ structures under the gauge-fixed Cartan
flow $(1.2)$. A closed $CR$ manifold $M$ is fillable if $M$ bounds a complex
manifold in the smooth $(C^{\infty})$ sense (i.e. there exists a complex
manifold with smooth boundary $M$, and the complex structure restricts to
the
$CR$ structure on $M$). The notion of fillability is weaker than that of
embeddability. Recall that a $CR$ manifold is embeddable if it can be
embedded
in $C^{N}$ for large $N$ with the $CR$ structure being the one induced from
the complex structure of $C^{N}$. The embeddability is a special property
for
3-dimensional $CR$ manifolds since any closed $CR$ manifold of dimension
$\geq$ 5 is embeddable ([BdM]). It is easy to see that a closed embeddable
(strongly pseudoconvex) $CR$ 3-manifold is fillable by some well-known
results
(see the argument on page 543 in [Ko]). Conversely, if there exists a smooth
strictly plurisubharmonic function defined in a neighborhood of a fillable
$M$, then $M$ is embeddable ([Ko], Theorem 5.3; in fact, any compact complex surface with nonempty strongly pseudoconvex boundary can be made Stein by deforming it and blowing down any exceptional curves according to [Bo]). Now it is natural to ask
the
following question:

\medskip

\paragraph{\textit{\ \ \ Is the embeddability (or fillability) preserved
under
the (gauge-fixed) Cartan flow (1.1) (or (1.2))?}}

\ \ \ \ \ \ \ \ \ \ \ \ \ \ \ \ \

\bigskip

An affirmative answer to the above question has an application in determining the topology of the space of all fillable $CR$ structures. For instance, one can apply such a result plus the convergence of the long time solution to (1.2) (expected for $S^3$) to prove that the space of all fillable $CR$ structures on $S^3$ is contractible (cf. Remark 4.3 in [El]). For other topological applications of solving (1.2), we refer the reader to [Ch].

Now by choosing a contact form, we can talk about the torsion of the
associated (positive) pseudohermitian structure (see \S2). Our first
observation is the following result.

\bigskip

\textit{Theorem A. Suppose there is a contact form }$\theta$\textit{ such
that
the torsion }$A_{J_{(0),}\theta}$\textit{ vanishes. Then under the
gauge-fixed
Cartan flow (1.2) (assuming smooth solution) with }$K=J_{(0)}$\textit{,
}$A_{J_{(t),}\theta}$\textit{ stays vanishing.}

\bigskip

Let $T$ denote the Reeb vector field associated with the contact form
$\theta$. The vanishing of the torsion is equivalent to saying that $T$ \ is
an infinitesimal $CR$ diffeomorphism (see (2.3) ). We say a $CR$ manifold
has
transverse symmetry if the infinitesimal generator of a one-parameter group
of
$CR$ diffeomorphisms is everywhere transverse to $\xi$. Such an
infinitesimal
generator can be realized as the Reeb vector field for a certain contact
form
$\hat{\theta}$ ([L2]). It follows from \textit{Theorem A} that

\bigskip

\textit{Corollary B. The }$CR$\textit{ structures }$J_{(t)}$\textit{ stay
having the same transverse symmetry as }$J_{(0)}$\textit{ does under the
gauge-fixed Cartan flow }$\mathit{(1.2)}$ \textit{with }$K=J_{(0)}$\textit{
and
}$\theta=\hat{\theta}.$

\bigskip

Let $W_{J,\theta}$ denote the Tanaka-Webster curvature of a
pseudohermitian structure $(J,\theta)$ (see \S2 for the definition). We have
the following result.

\bigskip

\textit{Theorem C. Suppose }$J_{(0)}$\textit{ is fillable with }%
$A_{J_{(0),}\theta}=0$\textit{ and }$W_{J_{(0)},\theta}>0$\textit{ (or }%
$<0$\textit{, respectively). Then the solution }$J_{(t)}$\textit{ to (1.2)
with }$K=J_{(0)}$\textit{ stays fillable for a short
time.}

\bigskip

Our proof of \textit{Theorem C }is a direct construction of an
integrable almost complex structure $\widetilde{J}$ on $M\times(0,\tau)$ for
a
small $\tau$ so that $\widetilde{J}|_{\xi}=J_{(t)}$ at $M\times\{t\}$ (see
\S 3 for details). Then we glue this complex structure $\widetilde{J}$ with
the one induced by the complex surface that $(M,J_{(0)})$ bounds along
$M\times\{0\}$ (identified with $M$). After we obtained the above result,
L\'{a}szl\'{o} Lempert pointed out to the author that the existence of a
$CR$
vector field $T$ is sufficient to imply the embeddability of the $CR$
structure. (see [Lem]) So by \textit{Theorem A} the condition\ in
\textit{Theorem C} on the Tanaka-Webster curvature can be removed according
to
[Lem]. We speculate that the embeddability (or fillability) is preserved
under
the (gauge-fixed) Cartan flow without any conditions.

On the other hand, the zero torsion condition reduces the complexity
of our flow a lot. It seems to be a good starting point. We are in a
situation
analogous to Hamilton's Ricci flow. Namely, given a closed contact
3-manifold
$(M,\xi)$, suppose there is a (positive) pseudohermitian structure
$(J,\theta)$ with vanishing torsion and positive Tanaka-Webster curvature.
Then can we deform $(J,\theta)$ according to the (gauge-fixed) Cartan flow
to
a limit $CR$ structure (together with the fixed contact form $\theta$) that
has the positive constant Tanaka-Webster curvature? It follows that this
limit
$CR$ structure has the vanishing Cartan tensor (recall that the torsion
stays
vanishing for all time). Therefore it is spherical.

In the Appendix (\S5), we show the uniqueness of the solution to (1.2)
with given smooth initial data.

The author would like to thank I-Hsun Tsai and L\'{a}szl\'{o} Lempert
for helpful conversations, and Jack Lee for e-mail communications. The proof
in
\S5 is based on the ideas of Jack, described in an e-mail message in January,
1995. This work was being done during the author's visit at the Institute
for
Advanced Study in the 2001-2002 academic year. He would therefore like to
thank the faculty and staff there for their hospitality during his stay.

\bigskip

\paragraph{2. Review in $CR$ and pseudohermitian geometry}

For most of basic material we refer the reader to [We], [Ta] or [L1].
Throughout the paper, our base space $M$ is a closed (compact with no
boundary) contact 3-manifold with the oriented contact structure $\xi$. A
$CR$
structure $J$ (compatible with $\xi$) is an endomorphism on $\xi$ with
$J^{2}=-identity$. By choosing a (global) contact form $\theta$ (exists if
the
normal bundle of $\xi$ in $TM$ is orientable), we can talk about
pseudohermitian geometry. The Reeb vector field $T$ is uniquely determined
by
${\theta}(T)=1$ and $T{\rfloor}d{\theta}=0.$ We choose a (local) complex
vector
field $Z_{1}$, an eigenvector of $J$ with eigenvalue $i$, and a (local)
complex 1-form ${\theta}^{1}$ such that $\{{\theta},{\theta^{1}},{\theta
^{\bar{1}}}\}$ is dual to $\{T,Z_{1},Z_{\bar{1}}\}$ (here
${\theta^{\bar{1}}}$
and $Z_{\bar{1}}$ mean the complex conjugates of ${\theta}^{1}$ and $Z_{1}$
respectively). It follows that $d{\theta}=ih_{1{\bar{1}}}{\theta^{1}}{\wedge
}{\theta^{\bar{1}}}$ for some nonzero real function $h_{1{\bar{1}}}$. If
$h_{1{\bar{1}}}$ is positive, we call such a pseudohermitian structure
$(J,{\theta})$ positive, and we can choose a $Z_{1}$ (hence $\theta^{1}$)
such
that $h_{1{\bar{1}}}=1$. That is to say

\bigskip

$(2.1)$ $\ \ \ \ \ \ d{\theta}=i{\theta^{1}}{\wedge}{\theta^{\bar{1}}}.$

\bigskip

We'll always assume our pseudohermitian structure $(J,{\theta})$ is
positive (by changing the sign of $\theta$ if negative) and
$h_{1{\bar{1}}}=1$
throughout the paper. The pseudohermitian connection of $(J,{\theta})$ is
the
connection ${\nabla}^{{\psi}.h.}$ on $TM{\otimes}C$ (and extended to
tensors)
given by%

$${\nabla}^{{\psi}.h.}Z_{1}={\omega_{1}}^{1}{\otimes}Z_{1},
{\nabla}^{{\psi}%
.h.}Z_{\bar1}={\omega_{\bar1}}^{\bar1}{\otimes}Z_{\bar1}, {\nabla}^{{\psi}%
.h.}T=0
$$

\noindent in which the connection 1-form ${\omega_{1}}^{1}$ is uniquely
determined by the following equation and associated normalization condition:

\bigskip

$(2.2)$ $\ \ \ \ \ \ \ \ \ \ \ \ \ \ \ \ d{\theta^{1}}={\theta^{1}}{\wedge
}{\omega_{1}}^{1}+{A^{1}}_{\bar{1}}{\theta}{\wedge}{\theta^{\bar{1}},}$

$\ \ \ \ \ \ \ \ \ \ \ \ \ \ \ \ \ \ \ \ \ \ \ {\omega_{1}}^{1}+{\omega
_{\bar{1}}}^{\bar{1}}=0.$

\bigskip

The coefficient ${A^{1}}_{\bar{1}}$ in (2.2) and its complex conjugate
${A^{\bar{1}}}_1$ are components of the torsion (tensor)
$A_{J,\theta}=i$
${A^{1}}_{\bar{1}}{\theta^{\bar{1}}\otimes}Z_{1}-i{A^{\bar{1}}}_{1}{\theta^{1}\otimes}Z_{\bar{1}}.$
Since $h_{1{\bar{1}}}=1$,
$A_{{\bar{1}}{\bar{1}}}=h_{1{\bar{1}}}{A^{1}}_{\bar{1}}={A^{1}}_{\bar{1}}$.
Further $A_{11}$ is just the complex conjugate of $A_{{\bar{1}}{\bar{1}}}%
$.\ Write $J=i$ ${\theta^{1}\otimes}Z_{1}-i{\theta^{\bar{1}}\otimes}Z_{\bar
{1}}$. It is not hard to see from (2.1) and (2.2) that

\bigskip

$(2.3)$\ \ \ $L_{T}J=2A_{J,\theta}$

\bigskip

\noindent where $L_{T}$ denotes the Lie differentiation in the direction $T$
(this is
the case when $f=-1$ in Lemma 3.5 of [CL1]). So the vanishing torsion is
equivalent to $T$ being an infinitesimal $CR$ diffeomorphism. We can define
the covariant differentiations with respect to the pseudohermitian
connection.
For instance, $f_{,1}=Z_{1}f$,
$f_{,1{\bar{1}}}=Z_{\bar{1}}Z_{1}f-{\omega_{1}%
}^{1}(Z_{\bar{1}})Z_{1}f$ for a (smooth) function $f$ (see, e.g., \S4 in
[L1]).
Now Differentiating ${\omega_{1}}^{1}$ gives

\bigskip

$(2.4)$ $\ \ \ \ \ \ \ \ \ \ \ \ \ \ d{\omega_{1}}^{1}=W{\theta^{1}%
}{\wedge}{\theta^{\bar{1}}}+2iIm(A_{11,{\bar{1}}}{\theta^{1}}{\wedge}{\theta})$

\bigskip

\noindent where $W$ or $W_{J,\theta}$ (to emphasize the dependence of the
pseudohermitian structure) is called the (scalar) Tanaka-Webster curvature.
There are distinguished $CR$ structures $J$, called spherical, if
$(M,\xi,J)$
is locally $CR$ equivalent to the standard 3-sphere $(S^{3},\hat{\xi
},\hat{J})$, or equivalently if there are contact coordinate maps into open
sets of $(S^{3},\hat{\xi})$ so that the transition contact maps can be
extended to holomorphic transformations of open sets in $C^{2}$. In 1930's,
Elie Cartan ([Ca], [CL1]) obtained a geometric quantity, denoted as $Q_{J}$,
by solving the local equivalence problem for $CR$ structure so that the
vanishing of $Q_{J}$ characterizes $J$ to be spherical. We will call $Q_{J}$
the Cartan (curvature) tensor. Note that $Q_{J}$ depends on a choice of
contact form $\theta.$ It is $CR$-covariant in the sense that if
$\tilde{\theta}=e^{2f}\theta$ is another contact form and $\tilde{Q}_{J}$ is
the corresponding Cartan tensor, then $\tilde{Q}_{J}=e^{-4f}Q_{J}$. We can
express $Q_{J}$ in terms of pseudohermitian invariants. Write $Q_{J}%
=iQ_{11}{\theta^{1}\otimes}Z_{\bar{1}}-iQ_{\bar{1}\bar{1}}{\theta^{\bar{1}%
}\otimes}Z_{1}$ (note that $Q_{1}{}^{\bar{1}}=Q_{11}$ and $Q_{\bar{1}}{}%
^{1}=Q_{\bar{1}\bar{1}}$ since we always assume $h_{1{\bar{1}}}=1$). We have
the following formula ([CL1], Lemma 2.2):

\bigskip

$(2.5)$\ \ \ \ \ $Q_{11}=\frac{1}{6}W_{,11}+\frac{i}{2}WA_{11}-A_{11,0}%
-\frac{2i}{3}A_{11,\bar{1}1.}$

\bigskip

In terms of local coframe fields we can express the Cartan flow (1.1) as
follows:

\bigskip

$(2.6)$ \ \ \ \ \ \ $\dot{\theta}^{1}=-Q_{\bar{1}\bar{1}}\theta^{\bar{1}}$

\bigskip

\noindent (cf. (2.16) in [CL1] with $E_{\bar{1}}$ $^{1}$ replaced by
$-iQ_{\bar{1}%
\bar{1}}$). The torsion evolves under the Cartan flow as shown in the follow
formula:

\bigskip

$(2.7)$ \ \ \ \ \ \ $\dot{A}_{11}=-Q_{11,0}$

\bigskip

\noindent (this is the complex conjugate of (2.18) in [CL1] with
$E_{\bar{1}}$ $^{1}$
replaced by $-iQ_{\bar{1}\bar{1}}$). Since the Cartan flow is invariant
under
the pullback action of contact diffeomorphisms (cf. the argument in the
proof
of Proposition 3.6 in [CL1]), we need to add a gauge-fixing term to the
right-hand side of (1.1) to get the subellipticity of its linearized
operator.
Let us recall what this term is. First we define a quadratic differential
operator $F_{J}$ from endomorphism fields to functions by ([CL1], p.236 and
note that $h_{1\bar{1}}=1$ here)

\bigskip

$(2.8)$\ \ \ \ \ \ \ \ $F_{J}E=(iE_{1\bar{1}}E_{11,\bar{1}\bar{1}}+iE_{\bar
{1}\bar{1}}E_{11,\bar{1}1})+conjugate.$

\bigskip

Also we define a linear differential operator $D_{J}$ from functions to
endomorphism fields and its formal adjoint $D_{J}^{\ast}$ by

\bigskip

$(2.9)$\ \ \ \ \ \ \ \ $D_{J}f=(f_{,11}+iA_{11}f){\theta^{1}\otimes}Z_{\bar
{1}}+(f_{,\bar{1}\bar{1}}-iA_{\bar{1}\bar{1}}f){\theta^{\bar{1}}\otimes}Z_{1},$

\ \ \ \ \ \ \ \ $\ \ \ \ \ \ \ D_{J}^{\ast}E=E_{11,\bar{1}\bar{1}}+E_{\bar
{1}\bar{1},11}-iA_{\bar{1}\bar{1}}E_{11}+iA_{11}E_{\bar{1}\bar{1}}.$

\bigskip

\noindent (note that we have used the notations $D_{J},D_{J}^{\ast}$ instead
of
$B_{J}^{^{\prime}},B_{J}$ in [CL1], respectively) Now let $K$ be a fixed
$CR$
structure. The Cartan flow with a gauge-fixing term reads as follows: (this
is (1.2))

\bigskip

\ \ \ \ \ \ \ \ \ \ $\frac{\partial J_{(t)}}{\partial
t}=2Q_{J_{(t)}}-\frac{1}%
{6}D_{J_{(t)}}F_{J_{(t)}}K.$

\bigskip

We also need the following commutation relations often:

\bigskip

$(2.10)$ \ \ \ \ \ \ \ $C_{I,01}-C_{I,10}=C_{I,\bar{1}}A_{11}-kC_{I}%
A_{11,\bar{1}}$

\ \ \ \ \ \ \ \ \ \ \ \ \ \ \ \ \ $C_{I,0\bar{1}}-C_{I,\bar{1}0}%
=C_{I,1}A_{\bar{1}\bar{1}}-mC_{I}A_{\bar{1}\bar{1},1}$

\ \ \ \ \ \ \ \ \ \ \ \ \ \ \ \
$C_{I,1\bar{1}}-C_{I,\bar{1}1}=iC_{I,0}+kC_{I}W.$

\bigskip

\noindent Here $C_{I}$ denotes a coefficient of a tensor with multi-index
$I$ consisting
of $1$ and $\bar{1},$ and $k$ is the number of $1$ in $I$ while $m$ is the
number of $\bar{1}$ in $I$ ([L2]).

\bigskip

\paragraph{3. Proof of Theorem A}

We'll compute the evolution of the torsion under the flow (1.2) (with
$K$ being the initial $CR$ structure $J_{(0)}$). First, instead of (2.7), we
have

\bigskip

$(3.1)$ \ \ \ \ \ $\dot{A}_{11}=-Q_{11,0}-\frac{i}{12}($
$D_{J}F_{J}K)_{11,0}.$

\bigskip

>From the formula (2.5) for $Q_{11}$, we compute $Q_{11,0}$. Using the
commutation relations (2.10) and the Bianchi identity: $W_{,0}=A_{11,\bar
{1}\bar{1}}+A_{\bar{1}\bar{1},11}$ ([L2]), we can express $Q_{11,0}$ only in
terms of $A_{11,}A_{\bar{1}\bar{1}}$ and their covariant derivatives as
follows:

\bigskip

$(3.2)$ \ \ \ \ $Q_{11,0}=\frac{1}{6}(A_{11,\bar{1}\bar{1}11}+A_{\bar{1}%
\bar{1},1111})-A_{11,00}-\frac{2i}{3}A_{11,\bar{1}10}+$ $l.w.t..$

\bigskip

\noindent where $l.w.t.$ means a lower weight term in $A_{11}$ and
$A_{\bar{1}\bar{1}}%
.$We count covariant derivatives in $1$ or $\bar{1}$ direction ($0$
direction,
resp.) as weight $1$ (weight $2$, resp.) and we call a term of weight $m$ if
its total weight of covariant derivatives is $m.$ For instance, $A_{11,\bar
{1}\bar{1}11}$, $A_{11,00}$ and $A_{11,\bar{1}10}$ are all of weight 4. So
more precisely each single term in $l.w.t.$ must contain terms of weight
$\leq3$ in $A_{11}$ or $A_{\bar{1}\bar{1}}.$In particular, if $A_{11}=0$,
then
$l.w.t.=0.$ Note that $A_{\bar{1}\bar{1},1111}$ is a ''bad''term in the
sense
that we need a gauge-fixing term to cancel it and obtain a fourth order
subelliptic operator in $A_{11}.$ Now by (2.9) the gauge-fixing term in
(3.1)
(up to a constant) reads as

\bigskip

$(3.3)$\ $\ \ \ \ \ ($ $D_{J}F_{J}K)_{11,0}=(F_{J}K)_{,110}+i[A_{11}%
(F_{J}K)]_{,0}$

$\ \ \ \ \ \ \ \ \ \ \ \ \ \ \ \ \ \ \ \ \ \ \ \ \ \ \ \ \ \ \ =$%
$(F_{J}K)_{,011}+l.w.t.$ (in $A_{11}$)

\bigskip

\noindent (we have used the commutation relations (2.10) for the last
equality). Write
$K=K_{11}\theta^{1}\otimes Z_{\bar{1}}+$ $K_{1\bar{1}}$\ ${\theta^{1}\otimes
}Z_{1}+K_{\bar{1}\bar{1}}$\ ${\theta^{\bar{1}}\otimes}Z_{1}+K_{\bar{1}%
1}{\theta^{\bar{1}}\otimes}Z_{\bar{1}}$ where
$K_{\bar{1}\bar{1}},K_{\bar{1}%
1}$ are the complex conjugates of $K_{11},K_{1\bar{1}}$, respectively. We
compute

\bigskip

$(3.4)$ \ \ \ \ \ $K_{11,0}=T[{\theta^{\bar{1}}(}KZ_{1})]-2{\omega_{1}}%
^{1}(T)K_{11}$

\ \ \ \ \ \ \ \ \ \ \ \ \ \ \ \ \ \ \ \  $=(L_{T}{\theta^{\bar{1}%
})(}KZ_{1})+{\theta^{\bar{1}}[(}L_{T}K)Z_{1}]+{\theta^{\bar{1}}[}K(L_{T}%
Z_{1})]-2{\omega_{1}}^{1}(T)K_{11}.$

\bigskip

It is easy to compute the first term using the (complex conjugate of)
structure equation (2.2) and the third term using the formula $[T,Z_{1}%
]=-A_{11}Z_{\bar{1}}+{\omega_{1}}^{1}(T)Z_{1}$ ([L1]). For the second term,
if
we take $K$ to be the initial $CR$ structure $J_{(0)},$then

\bigskip

$(3.5)$ \ \ \ \ \ $L_{T}K=L_{T}J_{(0)}=2A_{J_{(0)},\theta}=0$

\bigskip

\noindent by (2.3) and the assumption. So altogether we obtain

\bigskip

$(3.6)$ \ \ \ \ \ $K_{11,0}=A_{11}(K_{1\bar{1}}-K_{\bar{1}1})=2A_{11}%
K_{1\bar{1}}.$

\bigskip

Note that $K^{2}=-I$ implies that $K_{1\bar{1}}=\pm i(1+|K_{11}%
|^{2})^{\frac{1}{2}}$ and $K_{\bar{1}1}=-K_{1\bar{1}.}$ It follows that

\bigskip

$(3.7)$ \ \ \ \ \
$K_{1\bar{1},0}=-A_{11}K_{\bar{1}\bar{1}}+A_{\bar{1}\bar{1}%
}K_{11}.$

\bigskip

\noindent Here the point is that both $K_{11,0}$ and $K_{1\bar{1},0}$ are
linear in
$A_{11}$ and $A_{\bar{1}\bar{1}}$ with coefficients being ''0th-order'' in a
(co)frame. Using (3.6), (3.7), we can express $(F_{J}K)_{,011}$ as follows:

\bigskip

$(3.8)$ \ \ \ \ \ $(F_{J}K)_{,011}=iK_{1\bar{1}}K_{11,\bar{1}\bar{1}%
011}+iK_{\bar{1}\bar{1}}K_{11,\bar{1}1011}-iK_{\bar{1}1}K_{\bar{1}\bar
{1},11011}-$

\ \ \ \ \ \ \ \ \ \ \ \ \ \ \ \ \ \ \ \ \ \ \ \ \ \ \ \ \
$iK_{11}K_{\bar{1}\bar{1},1\bar{1}011}+l.w.t..$

\bigskip

\noindent Here and hereafter $l.w.t.$ will mean a lower weight term in
$A_{11}%
,A_{\bar{1}\bar{1}}$ up to weight 3 with coefficients in $K_{1\bar{1}},$
$K_{\bar{1}\bar{1}},$ $K_{\bar{1}1},$ $K_{11}$ and their covariant
derivatives
up to weight 5. Note that $A_{11},A_{\bar{1}\bar{1}}$ are of weight 2 in
$K_{1\bar{1}},K_{\bar{1}\bar{1}},$ $K_{\bar{1}1},$ $K_{11}.$ The first four
terms on the right-hand side of (3.8) contain the hightest weight terms of
weight 4 in $A_{11},A_{\bar{1}\bar{1}}$ in view of the commutation relations
(2.10) and (3.6), (3.7) as will be shown below. Using (2.10) repeatedly and
(3.6), we compute

\bigskip

$(3.9)$ \ \ \ \ \ $K_{11,\bar{1}\bar{1}011}=K_{11,\bar{1}0\bar{1}%
11}+l.w.t.=K_{11,0\bar{1}\bar{1}11}-2K_{11}A_{\bar{1}\bar{1},1\bar{1}11}+l.w.t.$

\ \ \ \ \ \ \ \ \ \ \ \ \ \ \ \ \ \ \ \ \ \ \ \   $=2K_{1\bar{1}}%
A_{11,\bar{1}\bar{1}11}-2K_{11}A_{\bar{1}\bar{1},1\bar{1}11}+l.w.t..$

\bigskip

Similarly we obtain

\bigskip

$(3.10)$ \ \ \ $K_{11,\bar{1}1011}$\ $=2K_{1\bar{1}}A_{11,\bar{1}111}%
-2K_{11}A_{\bar{1}\bar{1},1111}+l.w.t.$

\ \ \ \ \ \ \ \ \ \ \ \ $K_{\bar{1}\bar{1},11011}=2K_{\bar{1}1}A_{\bar
{1}\bar{1},1111}-2K_{\bar{1}\bar{1}}A_{11,\bar{1}111}+l.w.t.$

\ \ \ \ \ \ \ \ \ \ \ \  $K_{\bar{1}\bar{1},1\bar{1}011}=2K_{\bar{1}1}%
A_{\bar{1}\bar{1},1\bar{1}11}-2K_{\bar{1}\bar{1}}A_{11,\bar{1}\bar{1}11}+l.w.t..$

\bigskip

Substituting (3.9), (3.10) in (3.8), we get, in view of (3.3),

\bigskip

$\ (3.11)$ $\ \ ($ $D_{J}F_{J}K)_{11,0}=-2iA_{11,\bar{1}\bar{1}11}%
+2iA_{\bar{1}\bar{1},1111}+l.w.t..$

\bigskip

Now substituting (3.2) and (3.11) in (3.1) gives

\bigskip

$(3.12)$ \ \ \ $\dot{A}_{11}=-\frac{1}{3}A_{11,\bar{1}\bar{1}11}%
+A_{11,00}+\frac{2i}{3}A_{11,\bar{1}10}+l.w.t.$

\bigskip

\noindent (note that the ''bad''terms cancel). Define
$L_{\alpha}A_{11}=-A_{11,1\bar{1}%
}-A_{11,\bar{1}1}+i\alpha A_{11,0}$ for a complex number $\alpha$. Let
$L_{\alpha}^{\ast}$ be the formal adjoint of $L_{\alpha}.$ It is a direct
computation (cf. p1257 in [CL2]) that

\bigskip

$(3.13)$ \ \ \ \
$L_{\alpha}^{\ast}L_{\alpha}A_{11}=2(A_{11,11\bar{1}\bar{1}%
}+A_{11,\bar{1}\bar{1}11})-i(3+\alpha+\bar{\alpha}-|\alpha|^{2})A_{11,1\bar
{1}0}$

\ \ \ \ \ \ \ \ \ \ \ \ \ \ \ \ \ \ \ \ \ \ \ \ \ \ \ \ \ \ \ \ $+i(3-\alpha
-\bar{\alpha}-|\alpha|^{2})A_{11,\bar{1}10}+l.w.t..$

\bigskip

Using the commutation relations (2.10), we can easily obtain

\bigskip

$(3.14)$ \ \ \ \ $A_{11,11\bar{1}\bar{1}}=A_{11,\bar{1}\bar{1}11}%
+2iA_{11,\bar{1}10}+2iA_{11,1\bar{1}0}+l.w.t.$

\ \ \ \ \ \ \ \ \ \ \ \ \
$A_{11,00}=-iA_{11,1\bar{1}0}+iA_{11,\bar{1}10}+l.w.t..$

\bigskip

In view of (3.14) and (3.13), we can rewrite (3.12) as follows:

\bigskip

$(3.15)$ \ \ \ \
$\dot{A}_{11}=-\frac{1}{12}L_{\alpha}^{\ast}L_{\alpha}A_{11}+l.w.t.$

\bigskip

\noindent for $\alpha=4+i\sqrt{3}.$ Since $\alpha$ is not an odd integer,
$L_{\alpha}$
and hence $L_{\alpha}^{\ast}L_{\alpha}$ (note $L_{\alpha}^{\ast}%
=L_{\bar{\alpha}}$) are subelliptic (e.g. [CL1]). Taking the complex
conjugate
of (3.15) gives a similar equation for $A_{\bar{1}\bar{1}}$ only with
$\alpha$
replaced by $-\bar{\alpha}.$ On the other hand, we observe that $A_{11}%
=0,A_{\bar{1}\bar{1}}=0$ for all (valid) time is a solution to (3.15) and
its
conjugate equation (note that $l.w.t.$ vanishes if $A_{11}$ and $A_{\bar
{1}\bar{1}}$ vanish as remarked previously). Therefore by the uniqueness of
the solution to a (or system of) subparabolic equation(s) (in the Appendix
\S 5 , we give a proof of the uniqueness of the solution to (1.2). An
analogous argument works for a general subparabolic equation on a closed
manifold), we conclude that $A_{11}$ stays vanishing under the flow (1.2).

\bigskip

\paragraph{4. Proof of Theorem C}

Let $J_{(t)}$ be a solution to (1.2) for $0\leq t<\tau$ with given
initial $J_{(0)}$ having the property stated in Theorem C. We are going to
construct an almost complex structure $\check{J}$ \ on $M\times\lbrack
0,\tau),$ integrable on $M\times(0,\tau).$

There is a canonical choice of the (unitary) frame $Z_{1(t)}$ with
respect to $J_{(t)}$ ([CL1]). Write
$Z_{1(t)}=\frac{1}{2}(e_{1(t)}-ie_{2(t)})$
where $e_{1(t)},e_{2(t)}\in\xi$ and $J_{(t)}e_{1(t)}=e_{2(t)}.$ Let
$\{\theta,e_{(t)}^{1},e_{(t)}^{2}\}$ be a coframe dual to $\{T,e_{1(t)}%
,e_{2(t)}\}$ on $M.$ We'll identify $M\times\{t\}$ with $M$ (hence
$T(M\times\{t\})$ with $TM$). Now we define an almost complex structure
$\check{J}$ \ at each point in $M\times\{t\}$ as follows:

\bigskip

\ \ \ \ \ \ $\check{J}$\ $|_{\xi}=J_{(t)}$ , $\check{J}T=-a\frac{\partial
}{\partial t}+bT+a(\alpha e_{1(t)}+\beta e_{2(t)}).$

\bigskip

\noindent Here $a,b,\alpha,\beta$ are some real (smooth) functions of space
variable and
$t$, and $a\neq0$ (so $\check{J}$ $\frac{\partial}{\partial t}$ is
completely
determined from the above formulas and $\check{J}$ $^{2}=-identity.$
Strictly
speaking, $\alpha,\beta$ depend on the choice of frame while $a,b$ are
global). It is easy to see that the coframe dual to $\{e_{1(t)},$
$e_{2(t)},$
$\frac{\partial}{\partial t}-(b/a)T-\alpha e_{1(t)}-\beta e_{2(t)},$
$(1/a)T\}$ is $\{e_{(t)}^{1}+\alpha dt,$ $e_{(t)}^{2}+\beta dt,$ $dt,$
$a\theta+bdt\}.$ So the following complex 1-forms:

\bigskip

\ $(4.1)$\ \ \ \ \ $\Theta^{1}=(e_{(t)}^{1}+\alpha dt)+i(e_{(t)}^{2}+\beta
dt)=\theta_{(t)}^{1}+\gamma^{1}dt,$

\ $(4.2)$\ \ \ \ \ $\eta=(a\theta+bdt)-idt=a\theta+(b-i)dt$

\bigskip

\noindent are type (1,0) forms with respect to $\check{J}$ . Here
$\gamma^{1}%
=\alpha+i\beta$ is really the $Z_{1(t)}$ coefficient of the vector field
$\alpha e_{1(t)}+\beta e_{2(t)}.$ Let $\Lambda^{p,q}$ denote the space of
type
(p,q) forms. The integrability of $\check{J}$ is equivalent to $d\Lambda
^{1,0}\subset\Lambda^{2,0}+\Lambda^{1,1}$ or $\Lambda^{2,0}\wedge
d\Lambda^{1,0}=0.$ In terms of \ $\Theta^{1},\eta,$ the integrability
conditions read as follows:

\bigskip

$(4.3)$ \ \ \ \ \ $\eta\wedge$\ $\Theta^{1}\wedge d\eta=0,$

$(4.4)$ \ \ \ \ \ $\eta\wedge$\ $\Theta^{1}\wedge d\Theta^{1}=0.$

\bigskip

Substituting (4.1), (4.2) in (4.3) and making use of $d\theta
=d_{M}\theta=i\theta_{(t)}^{1}\wedge\theta_{(t)}^{\bar{1}}($here $d_{M}$
denotes the exterior differentiation on $M$ and $d=d_{M}+dt\frac{\partial
}{\partial t}$ on $M\times(0,\tau)),$ we obtain

\bigskip

$\ \ \ \ \ \ \ \ \ \ \ 0=\eta\wedge$\ $\Theta^{1}\wedge d\eta=[ab_{,\bar{1}%
}-(b-i)a_{,\bar{1}}+ia^{2}\gamma^{1}]\theta\wedge\theta_{(t)}^{1}\wedge
\theta_{(t)}^{\bar{1}}\wedge dt.$

\bigskip

\noindent Here $b_{,\bar{1}}=Z_{\bar{1}(t)}b,$
$a_{,\bar{1}}=Z_{\bar{1}(t)}a.$ Therefore
(4.3) is equivalent to the relation between $a,b$ and $\gamma^{1}$ as below:

\bigskip

$(4.5)$ \ \ \ $\gamma^{1}=ia^{-1}b_{,\bar{1}}-ia^{-2}(b-i)a_{,\bar{1}}.$

\bigskip

Next note that $d\theta_{(t)}^{1}=d_{M}\theta_{(t)}^{1}+dt\wedge
\dot{\theta}_{(t)}^{1}$ and

\bigskip

$\ \ \ \ \ \dot{\theta}_{(t)}^{1}=\{-Q_{\bar{1}\bar{1}(t)}+\frac{i}%
{12}[(F_{J_{(t)}}J_{(0)})_{,\bar{1}\bar{1}(t)}-iA_{\bar{1}\bar{1}%
(t)}(F_{J_{(t)}}J_{(0)})]\}\theta_{(t)}^{\bar{1}}.$

\bigskip

\noindent ($Q_{\bar{1}\bar{1}(t)}$ is the $\bar{1}\bar{1}-$component of the
Cartan
tensor with respect to $J_{(t)}$) So substituting (4.1),(4.2) in (4.4) and
making use of (2.2) for $\theta_{(t)}^{1}$, we obtain

\bigskip

\ \ \ \ $\ \ 0=\eta\wedge\Theta^{1}\wedge d\Theta^{1}=(a\{Q_{\bar{1}\bar
{1}(t)}-\frac{i}{12}[(F_{J_{(t)}}J_{(0)})_{,\bar{1}\bar{1}(t)}-iA_{\bar{1}%
\bar{1}(t)}(F_{J_{(t)}}J_{(0)})]\}$

$\ \ \ \ \ \ \ \ \ \ \ \ \ \ \ \ \ \ \ \ \ \ \ \ \ \ \ \ \ \ \ \ \
+(b-i)A_{\bar
{1}\bar{1}(t)}+a\gamma_{,\bar{1}}^{1})\theta\wedge\theta_{(t)}^{1}\wedge
\theta_{(t)}^{\bar{1}}\wedge dt.$

\bigskip

\noindent Here $A_{\bar{1}\bar{1}(t)}$ is the $\bar{1}\bar{1}-$component of
the
torsion tensor with respect to $J_{(t)}$ and $\gamma_{,\bar{1}}^{1}=Z_{\bar
{1}(t)}\gamma^{1}$ $+{\omega_{1}}_{(t)}^{1}(Z_{\bar{1}(t)})\gamma^{1}$ where
${\omega_{1}}_{(t)}^{1}$ is just the pseudohermitian connection form with
respect to $\theta_{(t)}^{1}$. Therefore (4.4) is equivalent to the
following
relation between $a,b$ and $\gamma_{,\bar{1}}^{1}:$

\bigskip

$(4.6)$ \ \ \ $\gamma_{,\bar{1}}^{1}=-Q_{\bar{1}\bar{1}(t)}+\frac{i}%
{12}[(F_{J_{(t)}}J_{(0)})_{,\bar{1}\bar{1}(t)}-iA_{\bar{1}\bar{1}%
(t)}(F_{J_{(t)}}J_{(0)})]-a^{-1}(b-i)A_{\bar{1}\bar{1}(t)}.$

\bigskip

Substituting (4.5) in (4.6) and letting $f=a^{-1},g=ba^{-1},u=f$ $+ig$,
we obtain an equation for a complex valued function $u:$

\bigskip

$(4.7)$ \ \ \ \ $u_{,\bar{1}\bar{1}}-iuA_{\bar{1}\bar{1}(t)}=-Q_{\bar{1}%
\bar{1}(t)}+\frac{i}{12}[(F_{J_{(t)}}J_{(0)})_{,\bar{1}\bar{1}(t)}-iA_{\bar
{1}\bar{1}(t)}(F_{J_{(t)}}J_{(0)})].$

\bigskip

In view of (2.9), we can express (4.7) in an intrinsic form:

\bigskip

$(4.8)$ \ \ \ \ $J_{(t)}\circ D_{J_{(t)}}f-D_{J_{(t)}}g=Q_{J_{(t)}}-\frac
{1}{12}D_{J_{(t)}}F_{J_{(t)}}J_{(0)}.$

\bigskip

Recall that $D_{J_{(t)}}f=\frac{1}{2}L_{X_{f}}J_{(t)}$ ([CL1]) in which
$X_{f}=-fT$ $+$ $i(Z_{1(t)}f)Z_{\bar{1}(t)}$ $-i(Z_{\bar{1}(t)}f)Z_{1(t)}$
is
the infinitesimal contact diffeomorphism induced by $f.$ So the image of
$D_{J_{(t)}}$ describes the tangent space of the orbit of the symmetry group
acting on $J_{(t)}$ by the pullback (in this case, the contact
diffeomorphisms
are our symmetries). Now (4.8) means that if $Q_{J_{(t)}}$ sits in the
''complexification'' of the infinitesimal symmetry group orbit for all
$t\in(0,\tau)$, then $\tilde{J}$ is integrable on $M\times(0,\tau).$

Now by \textit{Theorem A} we have $A_{11(t)}=0$ for $0\leq t<\tau.$
So in view of (2.5), we get

\bigskip

$(4.9)$ $\ \ \ \ \ Q_{11(t)}=\frac{1}{6}W_{,11(t)}.$

\bigskip

\noindent Here
$W_{,11(t)}=(Z_{1(t)}^{{}})^{2}W_{(t)}-{\omega_{1}}_{(t)}^{1}%
(Z_{1(t)})Z_{1(t)}W_{(t)}$ and $W_{(t)}$ is the Tanaka-Webster curvature
with
respect to $J_{(t)}$ (and fixed $\theta$). Therefore $a=-6W_{(t)}^{-1}$ and
$b=-\frac{1}{2}(F_{J_{(t)}}J_{(0)})W_{(t)}^{-1}$, hence $u=-\frac{1}{6}%
W_{(t)}+\frac{i}{12}F_{J_{(t)}}J_{(0)}$ is a solution to (4.7) by (4.9) for
$0\leq t<\tau$ with $\tau$ small so that $W_{(t)}>0$ or $W_{(t)}<0.$ Thus
for
such a choice of $a,b,$ and $\tau,$ our $\check{J}$ is integrable on
$M\times(0,\tau).$

On the other hand, $(M,J_{(0)})$ bounds a complex surface $N$ by our
assumption. So we have another almost complex structure $\hat{J}$ on
$M\times(-\delta,0]$ induced from $N$, integrable on $M\times(-\delta,0),$
and
restricting to $J_{(0)}$ on $(M,\xi).$Up to a diffeomorphism from
$M\times(-\delta_{1},0]$ to $M\times(-\delta_{2},0],$ identity on
$M\times\{0\}$ for $\delta_{1,}\delta_{2}$ perhaps smaller than $\delta,$ we
can assume that $\check{J}$ and $\hat{J}$ coincide at $M\times\{0\}$ where
they may not coincide up to $C^{k}$ for $k\geq1,$ however. We want to find a
local diffeomorphism $\Phi$ from a neighborhood $U$ of a point in $M$ times
$(-\delta_{1},0]$ to a similar set so that $\Phi$ is an identity on
$U\times\{0\},$ and $\Phi^{\ast}\hat{J}$ coincides with $\check{J}$ up to
$C^{k}$ for some large integer $k$ at $U\times\{0\}$. Let $x^{i},0\leq
i\leq3$
denote the coordinates of $U\times(-\delta_{1},0]$ with $x^{0}$ being the
time
variable for $(-\delta_{1},0].$ Let $y^{i},0\leq i\leq3$ denote the
corresponding coordinates of the image of $\Phi$ with $y^{0}$ being the time
variable. If we express $\hat{J}_{1}=\Phi^{\ast}\hat{J}$ $=\Phi_{\ast}%
^{-1}(\hat{J}\circ\Phi)\Phi_{\ast}$\ in coordinates, we usually write

\bigskip

\ \ $(\hat{J}_{1})_{m}^{l}=\hat{J}_{i}^{j}\frac{\partial y^{i}}{\partial
x^{m}}\frac{\partial x^{l}}{\partial y^{j}}$

\bigskip

\noindent for $\hat{J}_{1}=$\
$(\hat{J}_{1})_{m}^{l}dx^{m}\otimes\frac{\partial
}{\partial x^{l}}$ and $\hat{J}=$
$\hat{J}_{i}^{j}dy^{i}\otimes\frac{\partial
}{\partial y^{j}}.$ Let $\eta=\Phi_{\ast}^{-1}.$ Then $\eta^{-1}$ has the
expression $(\frac{\partial y^{i}}{\partial x^{m}}),$ the Jacobian matrix of
$\Phi,$ in coordinates. We require $\eta^{-1}=identity$ at each point with
$x^{0}=0$ where $\hat{J}$ coincides with $\check{J}$. Differentiating
$\hat{J}_{1}=\Phi^{\ast}\hat{J}$ $=\Phi_{\ast}^{-1}(\hat{J}\circ\Phi
)\Phi_{\ast}=\eta(\hat{J}\circ\Phi)\eta^{-1}$(considered as a matrix
equation
with respect to the above-mentioned bases) in $x^{0}$ at $x^{0}=0$, we
obtain

\bigskip

$(4.10)$ $\ \ \ \hat{J}_{1}^{\prime}-\hat{J}^{\prime}=\eta^{\prime
}\hat{J}-\hat{J}\eta^{\prime}.$

\bigskip

\noindent Here the prime of $\hat{J}^{\prime}$ means the $y^{0}$-derivative
at $y^{0}=0$
while the prime of $\hat{J}_{1}^{\prime}$ and $\eta^{\prime}$means the
$x^{0}-$derivative at $x^{0}=0$. Finding $\Phi$ such that $\hat{J}_{1}%
=\Phi^{\ast}\hat{J}$ \ coincides with $\check{J}$ up to $C^{1}$ at
$U\times\{0\}$ is reduced to solving the above equation $(4.10)$ for
$\eta^{\prime}$ with $\hat{J}_{1}^{\prime}=\check{J}^{\prime}.$ Here the
prime
of $\check{J}^{\prime}$ means the $t$-derivative at $t=0.$ And this can be
done by simple linear algebra as follows. First note that $C=\check{J}%
^{\prime}-\hat{J}^{\prime}$ satisfies $\hat{J}C+C\hat{J}=0$ since $\check
{J}=\hat{J}$ at $U\times\{0\}$ and both $\check{J}^{\prime}$ and
$\hat{J}^{\prime}$ satisfies the same relation as $C$ does. With respect to
a
suitable basis, $\hat{J}$ has a canonical matrix representation:

$\left(
\begin{array}
[c]{cccc}%
0 & -1 & 0 & 0\\
1 & 0 & 0 & 0\\
0 & 0 & 0 & -1\\
0 & 0 & 1 & 0
\end{array}
\right)  .$

\noindent Then $C$ has the matrix form $\left(
\begin{array}
[c]{cc}%
C_{11} & C_{12}\\
C_{21} & C_{22}%
\end{array}
\right)  $ where each $C_{ij}$ is a $2\times2$ matrix $\left(
\begin{array}
[c]{cc}%
a_{ij} & b_{ij}\\
b_{ij} & -a_{ij}%
\end{array}
\right)  .$ Now the solution $\eta^{\prime}$ to $(4.10)$ has the matrix form
$\left(
\begin{array}
[c]{cc}%
\eta_{11}^{\prime} & \eta_{12}^{\prime}\\
\eta_{21}^{\prime} & \eta_{22}^{\prime}%
\end{array}
\right)  $ where each $\eta_{ij}^{\prime}$ is a $2\times2$ matrix $\left(
\begin{array}
[c]{cc}%
u_{ij} & v_{ij}\\
w_{ij} & s_{ij}%
\end{array}
\right)  $ satisfying the relations: $v_{ij}+w_{ij}=-a_{ij},$ $u_{ij}%
-s_{ij}=b_{ij}.$ Once $\eta^{\prime}$ is determined by the equation (4.10),
it
is easy to construct the ''local''diffeomorphism $\Phi_{1}$ such that the
inverse Jacobian and its $x^{0}$-derivative at $x^{0}=0$ of $\Phi_{1}$ is
$\eta=$the identity and $\eta^{\prime}$, respectively (we may need to shrink
the time interval $(-\delta_{1},0]$). So if we start with $\hat{J}_{1}%
=\Phi_{1}^{\ast}\hat{J}$ instead of $\hat{J}$ and repeat the above procedure
looking for $\Phi_{2}$ so that $\hat{J}_{2}=\Phi_{2}^{\ast}\hat{J}_{1}$
coincides with $\check{J}$ at $U\times\{0\}$ up to $C^{2}$, we differentiate
$\hat{J}_{2}=\eta_{1}\hat{J}_{1}\eta_{1}^{-1}$ twice with respect to $x^{0}%
$\bigskip\ at $x^{0}=0.$ Here $\eta_{1}$ denotes the inverse Jacobian matrix 
of
$\Phi_{2}$ (to be determined). Requiring $\hat{J}_{2}^{\prime}%
=\hat{J}_{1}^{\prime}$ and $\eta_{1}=identity$ (at $x^{0}=0$) implies
$\eta_{1}^{\prime}=0.$ It then follows that $\eta_{1}^{\prime\prime},$the
second derivative of $\eta_{1}$ in $x^{0}$ at $x^{0}=0,$ satisfies a similar
equation as in (4.10):

\bigskip

$(4.11)$ \ \ \ $\eta_{1}^{\prime\prime}\hat{J}_{1}-\hat{J}_{1}\eta_{1}%
^{\prime\prime}=\hat{J}_{2}^{\prime\prime}-\hat{J}_{1}^{\prime\prime}.$

\bigskip

Now we can verify that the right-hand side anti-commutes with $\hat{J}_{1}$
as
follows:
$(\hat{J}_{2}^{\prime\prime}-\hat{J}_{1}^{\prime\prime})\hat{J}_{1}%
+$\ $\hat{J}_{1}(\hat{J}_{2}^{\prime\prime}-\hat{J}_{1}^{\prime\prime
})=(\hat{J}_{2}^{\prime\prime}\hat{J}_{1}+\hat{J}_{1}\hat{J}_{2}^{\prime
\prime})-(\hat{J}_{1}^{\prime\prime}\hat{J}_{1}+\hat{J}_{1}\hat{J}_{1}%
^{\prime\prime})=-2(\hat{J}_{2}^{\prime})^{2}+2(\hat{J}_{1}^{\prime})^{2}=0$
(here we have used $\hat{J}_{2}=\hat{J}_{1}$, $\hat{J}_{2}^{\prime
}=\hat{J}_{1}^{\prime}$ and $J^{\prime\prime}J+2(J^{\prime})^{2}%
+JJ^{\prime\prime}=0$ for any almost complex structure $J$ by
differentiating
$J^{2}=-I$ twice. So we can solve (4.11) for $\eta_{1}^{\prime\prime}$ with
$\hat{J}_{2}^{\prime\prime}=\check{J}^{\prime\prime}$ and hence find a
$\Phi_{2}$ with the required properties as before. In general, suppose we
have
found $\Phi_{n-1}$such that $\hat{J}_{n-1}=$
$\Phi_{n-1}^{\ast}\hat{J}_{n-2}%
=\eta_{n-2}\hat{J}_{n-2}\eta_{n-2}^{-1}$ coincides with $\check{J}$ up to
$C^{n-1}$ at $x^{0}=0.$ Then by the similar procedure we can find $\Phi_{n}$
such that $\hat{J}_{n}=\Phi_{n}^{\ast}\hat{J}_{n-1}=\eta_{n-1}\hat{J}_{n-1}%
\eta_{n-1}^{-1}$ coincides with $\check{J}$ up to $C^{n}$ at $x^{0}=0,$ and
the $x^{0}$-derivatives of $\eta_{n-1}$ vanish up to the order $n-1.$
Furthermore the $n$-th $x^{0}$-derivative $\eta_{n-1}^{(n)}$ satisfies a
similar equation as in (4.10) or (4.11):

\bigskip

$(4.12)$ \ \ \ $\eta_{n-1}^{(n)}\hat{J}_{n-1}-\hat{J}_{n-1}\eta_{n-1}%
^{(n)}=\check{J}^{(n)}-\hat{J}_{n-1}^{(n)}.$

\bigskip

\noindent Here $\check{J}^{(n)}$ denotes the $n$-th $t$-derivative of
$\check{J}$ at
$t=0$ while $\hat{J}_{n-1}^{(n)}$ means the $n$-th $x^{0}$-derivative of
$\hat{J}_{n-1}$ at $x^{0}=0$ $.$

Now $\hat{J}_{n}$ defined on $U\times(-\delta_{n},0]$ and $\check{J}$
\ defined on $U\times\lbrack0,\delta_{n})$ for a small $\delta_{n}>0$
together
form a $C^{n}$ integrable almost complex structure on $U\times(-$ $\delta
_{n},\delta_{n})$. Therefore $U\times(-$ $\delta_{n},\delta_{n})$ is a
complex
manifold for $n\geq4$ by a theorem of Newlander-Nirenberg ([NN]). Since $M$
is
compact, we can cover it by a finite number of $U^{\prime}s$ and have
corresponding $\delta_{n}^{\prime}s$. For each point in the overlap of two
$U^{\prime}s$ considered in $U\times\{0\}$, we can find local coordinate
maps
from an open neighborhood $V$ contained in the intersection of two
associated
$U\times(-$ $\delta_{n},\delta_{n})^{\prime}s$ into $C^{2}$ so that the
transition map $\psi$ on the ''concave'' part corresponding to positive
''time
variable'' is holomorphic (note that our $(M,J_{(0)})$ is a strongly
pseudoconvex boundary of $N.$ And on $V\cap\{M\times\lbrack0,\tau)\},$ we
have
the ''same'' integrable almost complex structure $\check{J}$ while on the
intersection of $V$ and $U\times(-$ $\delta_{n},0)^{\prime}s,$ we may have
''different'' $(\hat{J}_{n})^{\prime}s$). We then extend $\psi$ to the
pseudoconvex part holomorphically, and denote the extension map by
$\tilde{\psi}$. Now glue $V\cap\{U\times(-$ $\delta_{n},0)\}$ (complex
structure $\hat{J}_{n}$) with $V\cap\{$another copy of $U\times(-$ $\delta
_{n},0)\}$ (perhaps different $\hat{J}_{n}$) through $\tilde{\psi}.$ In this
way we can manage to extend the complex structure $\check{J}$ across
$M\times\{0\}$ to $M\times(-\delta,0)$ (globally) for some small $\delta>0$.
Finally the identity (a $CR$ diffeomorphism) on $(M,J_{(0)})$ extends to a
biholomorphism $\rho$ between $M\times(-\delta,0)$ (perhaps smaller
$\delta$)
and an open set in $N$ near $M$ (recall that $N$ is a complex surface that
$M$
bounds). Glue $M\times(-\delta,t)$ $(t<\tau)$ and $N$ \ via $\rho$ to form a
complex surface $N_{t}$ that ($M\times\{t\},J_{(t)})$ bounds. We have shown
that $J_{(t)}$ is fillable for $0<t<\tau.$

\paragraph{5. Appendix: uniqueness of the solution to (1.2)}

In this section, we'll show that the short-time solution ([CL1]) to the
gauge-fixed Cartan flow $(1.2)$ with given initial data is actually unique.
As mentioned in Introduction, the idea of proof was suggested by Jack Lee.

First we refer the reader to [CL1] for the definitions of various
notations, e.g., the Folland-Stein space $S_{k}$, some time-dependent space
$E_{k,\tau},$ the vector bundle ${\cal{E}}_{J},$ the operators
$L_{\alpha},\Lambda$,
etc.. We define the space $\tilde{A}_{k+4,\tau}$ to consist of all the
elements $u$ in $E_{k+4,\tau}$ with the initial value $u_{(0)}$ in
$S_{k+2}.$
Let $\Xi_{k,\tau}=\{(u_{(0)},u)\in S_{k+2}\times\tilde{A}_{k+4,\tau}$ $|$
$u\in\tilde{A}_{k+4,\tau}\}.$ Let $P_{(t)}$ be a time-dependent linear
operator on sections of ${\cal{E}}_{J}$ over $M,$ involving only spatial
derivatives
of weight $\leq4$ and depending smoothly on $t\in\lbrack0,1],$ such that
$P_{(0)}=cL_{\alpha}^{\ast}L_{\alpha}+S$, where $c$ is a positive constant,
$\alpha$ is admissible and $S$ is an operator of weight $\leq3$ (here
$L_{\alpha}^{\ast}$ instead of $L_{\alpha}^{^{\prime}}$ in [CL1] means the
adjoint operator of $L_{\alpha}$). The following theorem extends Theorem 4.6
in [CL1] to the case of nonvanishing initial data $u_{(0)}.$

\bigskip

\textit{Theorem 5.1. For any integer }$k\geq0,$\textit{ there exists
}$0<\tau\leq1$\textit{ such that the map }$\Psi:\Xi_{2k,\tau}\rightarrow
S_{2k+2}\times E_{2k,\tau}$\textit{ defined by }

\bigskip

\textit{ \ \ \ }$\Psi((u_{(0)},u))=(u_{(0)},(\partial_{t}+P_{(t)})u)$

\bigskip

\noindent \textit{is a bounded isomorphism.}

\bigskip

$Proof.$ \ \ We will follow the treatment given in the proof of
Theorem 4.6 in [CL1]. It is clear that $\Psi$ is linear and bounded. $\Psi$
being surjective is equivalent to solving the following initial-value
problem:

\bigskip

$(5.1)$ \ \ \ \ \ \ \ \ \ $(\partial_{t}+P_{(t)})u_{(t)}=f_{(t)},$

$(5.2)$ \ \ \ \ \ \ \ \ \ \ \ \ \ \ \ \ \ \ $u_{(0)}=g$

\bigskip

\noindent for $u\in\tilde{A}_{2k+4,\tau}$ with given $(g,f)\in
S_{2k+2}\times
E_{2k,\tau}.$ Recall ([CL1], p.244) that $B_{k,\tau}=\{u\in E_{k,\tau}%
:u_{(\tau)}=0\}.$ Define $\tilde{B}_{k,\tau}=\{u\in B_{k,\tau}:u_{(0)}\in
S_{k-2}\}.$ Let $\Sigma_{k,\tau}=\{(v_{(0)},v)\in S_{k-2}\times\tilde
{B}_{k,\tau}:v\in\tilde{B}_{k,\tau}\}.$ Define a hermitian bilinear form
$\Omega:\Sigma_{2k+4,\tau}\times(S_{2k+2}\times S_{2k+4,\tau})\rightarrow C$
by

\bigskip

$(5.3)$ \ \ $\Omega((v_{(0)},v),(h,u))=A(v,u).$

\bigskip

\noindent Here $A(v,u)$ is just the hermitian bilinear form that we used in
[CL1] (see
(4.9) on page 245). Observe that

\ \ \ (1) For any $(v_{(0)},v)\in\Sigma_{2k+4,\tau},$ the linear functional
$(h,u)\rightarrow\overline{\Omega((v_{(0)},v),(h,u))}$ is bounded on
$S_{2k+2}\times S_{2k+4,\tau}$ since $|\overline{A(v,u)}|\leq
C||u||_{2k+4,\tau}\leq C(||h||_{2k+2}+||u||_{2k+4,\tau}).$

\ \ \ (2) For some positive constant $C,$ $C(||v_{(0)}||_{2k+2}^{2}%
+||v||_{2k+4,\tau}^{2})\leq|\Omega((v_{(0)},v),(v_{(0)},v))|$ for all
$(v_{(0)},v)\in\Sigma_{2k+4,\tau}.$ To verify this, we review the argument
at the bottom of page 245 and the top of page 246 in [CL1] and conclude that

\begin{eqnarray}
|\Omega((v_{(0)},v),(v_{(0)},v))|&=&|A(v,v)|{\geq}Re{\:} A(v,v)\nonumber\\
&\geq
&C^{\prime}||v||_{2k+4,\tau}^{2}+\frac{1}{2}||\Lambda^{k+1}v_{(0)}||_{0}^{2}\nonumber\\
&\geq
&C(||v||_{2k+4,\tau}^{2}+||v_{(0)}||_{2k+2}^{2}).\nonumber
\end{eqnarray}

\noindent Here $C^{\prime}$ and $C$ are some positive constants, and the
last inequality
follows from Corollary 4.3 in [CL1]. Under conditions (1),(2), we can apply
a
generalized Lax-Milgram lemma due to J. L. Lions ([Tr], lemma 41.2) to
assert
that for any continuous linear functional $G$ on $S_{2k+2}\times
S_{2k+4,\tau
},$ there exists $(\tilde{h},\tilde{u})\in S_{2k+2}\times S_{2k+4,\tau}$
with

\bigskip

$(5.4)$ \ \ \ \ \ $||\tilde{h}||_{2k+2}+||\tilde{u}||_{2k+4,\tau}\leq
C||G||$
\ \ \ (operator norm)

\bigskip

\noindent such that

\bigskip

$(5.5)$ \ \ \ \ \ $\Omega((v_{(0)},v),(\tilde{h},\tilde{u}))=G((v_{(0)},v))$
\ \ \ for all $(v_{(0)},v)\in\Sigma_{2k+4,\tau}.$

\bigskip

Now given $(g,f)\in S_{2k+2}\times E_{2k,\tau},$ we take $G:$
$S_{2k+2}\times S_{2k+4,\tau}\rightarrow C$ to be the functional

\bigskip

$(5.6)$\ \ \ $G((h,v))=F(v)+(\Lambda^{k+1}h,\Lambda^{k+1}g)_{J}$

\bigskip

\noindent in which $F:S_{2k+4,\tau}\rightarrow C$ is given by ([CL1], p.246)

\bigskip

$F(v)=\int_{0}^{\tau}(\Lambda^{k+2}v_{(t)},\Lambda^{k}(e^{-\kappa
t}f_{(t)}))_{J}dt.$

\bigskip

It is easy to see that $|G((h,v))|\leq
C_{1}(||h||_{2k+2}+||v||_{2k+4,\tau})$
since $|F(v)|\leq C_{2}||v||_{2k+4,\tau}.$ Here $C_{1},C_{2}$ are some
positive constants. Thus there exists $(\tilde{h},\tilde{u})\in$
$S_{2k+2}\times S_{2k+4,\tau}$ satisfying $(5.4)$ so that $(5.5)$ holds. By
taking $v\in C_{0}^{\infty}((0,\tau)\times M)$ (smooth with compact
support)(which implies $v_{(0)}=0$), we are reducing the equation $(5.5)$ to
the equation $(4.12)$ in [CL1]. So an argument there on p.246 shows that
$u_{(t)}=e^{\kappa t}\tilde{u}_{(t)}$ satisfies $(5.1).$ Furthermore, we can
show that $u\in E_{2k+4,\tau}$ by $(5.4)$ and a similar argument as in
[CL1],
p.247. Now for $v\in B_{2k+4,\tau},$ we have

\bigskip

$(5.7)\ \ \ \ \ \ \ \ \ \ \ \ \ \ \ 0=A(v,\tilde{u})-G((v_{(0)},v))$
\ \ \ \ \ \ (by $(5.5),(5.3)$)

\ \ \ \ \ \ \ \ \ \ \ \ \ \ \ \ \ \ \ \ \ \ \ \
$=(A(v,\tilde{u})-F(v))-(\Lambda
^{k+1}v_{(0)},\Lambda^{k+1}g)_{J}$ \ \ (by $(5.6)$)

\ \ \ \ \ \ \ \ \ \ \ \ \ \ \ \ \ \ \ \ \ \ \ \  $=(\Lambda^{k+1}%
v_{(0)},\Lambda^{k+1}\tilde{u}_{(0)})_{J}-(\Lambda^{k+1}v_{(0)},\Lambda
^{k+1}g)_{J}$

\ \ \ \ \ \ \ \ \ \ \ \ \ \ \ \ \ \ \ \ \ \ (see the formula in the middle
of
p.247 in [CL1])

\bigskip

Since, for any $w\in C^{\infty}(M),$ we can find $v\in B_{2k+4,\tau}$ such
that $v_{(0)}=w,$ we obtain from $(5.7)$ that $\Lambda^{2k+2}(\tilde{u}%
_{(0)}-g)=0$ in the distribution sense, and hence
$u_{(0)}=\tilde{u}_{(0)}=g.$
We have shown that $u\in\tilde{A}_{2k+4,\tau}$ and is a solution to
$(5.1)-(5.2).$ Therefore $\Psi$ is surjective. By the same proof for the
uniqueness as in [CL1], p.247 (the last paragraph of the proof for Theorem
4.6), we conclude that $\Psi$ is injective.

\ \ \ \ \ \ \ \ \ \ \ \ \ \ \ \ \ \ \ \ \ \ \ \ \ \ \ \ \ \ \ \ \ \ \ \ \ \
\ \ \ \ \ \ \ \ \ \ \ \ \ \ \ \ \ \ \ \ \ \ \ \ \ \ \ \ \ \ \ \ \ \ \ \ \ \
\ \ \ \ \ \ \ \ \ \ \ \ \ \ Q.E.D.
\ \ \

\bigskip

We remark that Theorem 5.1 may still be valid for any positive
half-integer $k$ if we can make sense of $\Lambda^{k}$, say, in terms of
pseudodifferential operators. In the elliptic case, the (complex) power of
an
elliptic operator (acting on sections of a vector bundle over a closed
manifold) can be well defined (see, e.g., [See]). For our case, one expects
to
define the power of a subelliptic operator like $\Lambda$ along the same
line
of ideas as in [See] together with the symbol calculus for so-called
$V$-operators ([BG]). On the other hand, Theorem 5.1 is sufficient for our
purpose of proving the uniqueness of the solution to (1.2)$.$

\bigskip

\textit{Theorem 5.2. Let }$\hat{J}$\textit{ be any smooth
}$(i.e. C^{\infty})$\textit{ oriented }$CR$\textit{ structure on
}$M.$\textit{
For a large enough integer }$m,$\textit{ say }$m\geq14,$\textit{ suppose
}$J_{(t)}^{1},J_{(t)}^{2}$\textit{ are two }$C^{m}$\textit{ solutions to
}$(1.2)$\textit{ on some time interval }$[0,\varepsilon]$\textit{ with
}$J_{(0)}^{1}=J_{(0)}^{2}=\hat{J}$\textit{. Then }$J_{(t)}^{1}=J_{(t)}^{2}%
$\textit{ on }$[0,d]$\textit{ for a small positive
}$d<\varepsilon$\textit{.}

\bigskip

$Proof.$ \ \ Recall ([CL1]) that we can parametrize $J_{(t)}$
by a section $u_{(t)}$ of ${\cal{E}}_{\hat{J}}$ and describe the nonlinear
operator
$\mathit{P}(J)=-2Q_{J}+\frac{1}{6}D_{J}F_{J}K$ by $\mathit{\hat{P}}$ from
sections of ${\cal{E}}_{\hat{J}}$ to sections of ${\cal{E}}_{\hat{J}}.$ The
equation
$(1.2)$ together with the initial condition $J_{(0)}=\hat{J}$ is equivalent
to

\bigskip

$(5.8)$ \ \ \ \ $\partial_{t}u_{(t)}+\mathit{\hat{P}}(u_{(t)})=0,$

$(5.9)$ \ \ \ \ \ \ \ \ \ \ \ \ \ $u_{\left(  0\right)  }=0$

\bigskip

\noindent ([CL1], (5.1),(5.2)). Let $\tilde{u}_{(t)}$ be an infinite-order
formal
solution to $(5.8)-(5.9)$ so that $\tilde{f}=(\partial_{t}+$
$\mathit{\hat{P}%
})\tilde{u}$ vanishes to infinite order at $t=0.$ Let $P$ denote the
linearization of $\mathit{\hat{P}}:\tilde{A}_{2k+4,\tau}\rightarrow
E_{2k,\tau}$ about $\tilde{u}.$ Observe that $P_{(0)}$ satisfies the
required
property, and $\mathit{\hat{P}}$ is $C^{1}$ for $k$ large enough, say,
$k\geq5$ ([CL1], pp.250-251). So in view of \textit{Theorem 5.1,} we can 
apply
the
inverse function theorem to conclude that

\bigskip

\textit{Lemma 5.3. For }$k\geq5,$\textit{ there exists
}$0<\tau\leq1$\textit{
such that the map }$Id\times(\partial_{t}+\hat{P}):\Xi_{2k,\tau}\rightarrow
S_{2k+2}\times E_{2k,\tau}$\textit{ defined by}

\bigskip

\textit{ \ \ \ \ \ \
}$(u_{(0)},u)\rightarrow(u_{(0)},(\partial_{t}+\hat{P})u)$

\bigskip

\noindent \textit{has a }$C^{1}$\textit{ inverse on some neighborhood of }%
$(0,\tilde{u})$\textit{ in }$\Xi_{2k,\tau}.$

\bigskip

Note that the smooth section $f_{\varepsilon}$ given by
\ $f_{\varepsilon(t)}=0$ for $0\leq t\leq\varepsilon,$ and $f_{\varepsilon
(t)}=\tilde{f}_{(t-\varepsilon)}$ for $\varepsilon\leq t\leq\tau$ is
arbitrarily close to $\tilde{f}$ in $E_{2k,\tau}$ for small $\varepsilon>0.$
So by Lemma 5.3 there exists $u\in A_{2k+4,\tau}$ $(u_{(0)}=0)$ satisfying
$(\partial_{t}+\mathit{\hat{P}})u=f_{\varepsilon}.$ For $t\in\lbrack
0,\varepsilon],$ u is a solution to $(5.8)-(5.9).$

Now suppose $v$ $\in A_{2k+4,\tau}$ is another solution to
$(5.8)-(5.9)$ also for $t\in\lbrack0,\varepsilon].$ Observe that
$f_{\varepsilon-d}$ is arbitrarily close to $f_{\varepsilon}$ in
$E_{2k,\tau}$
for small $d>0.$ Moreover, $(v_{(d)},f_{\varepsilon-d})$ is in the
neighborhood of $(0,\tilde{f})$ where we can apply Lemma 5.3 if $d$ is small
enough. So by Lemma 5.3 there exists $\tilde{v}$ close to $\tilde{u}$ in
$\tilde{A}_{2k+4,\tau}$ such that

\bigskip

$(5.10)$ \ \ \ \ \
$(\partial_{t}+\mathit{\hat{P}})\tilde{v}=f_{\varepsilon-d},$

$(5.11)$ \ \ \ \ \ \ \ \ \ \ \ \ \ $\tilde{v}_{(0)}=v_{(d)}.$

\bigskip

Let $w_{(t)}=v_{(t)}$ for $0\leq t\leq d$ and $w_{(t)}=\tilde
{v}_{(t-d)}$ for $d\leq t\leq\tau.$ The equation insures that the
derivatives of $v$ and $\tilde{v}$ match up at $t=d$. Compute $(\partial
_{t}+\mathit{\hat{P}})w=0$ for $0\leq t\leq d$ and $(\partial_{t}%
+\mathit{\hat{P}})w=f_{\varepsilon}$ for $d\leq t\leq\tau.$ (note that
$f_{\varepsilon-d(t-d)}=f_{\varepsilon(t)}$) Since $d<\varepsilon,$ we
actually have $(\partial_{t}+\mathit{\hat{P}})w=f_{\varepsilon}.$ Also, for
$d$ small, $w$ is close to $\tilde{u}$ since $v_{(0)}=0$ and $\tilde{v}$ is
close to $\tilde{u}.$ Therefore $w=u.$ It follows that $v_{(t)}=u_{(t)}$ for
$t\in\lbrack0,d].$

\ \ \ \ \ \ \ \ \ \ \ \ \ \ \ \ \ \ \ \ \ \ \ \ \ \ \ \ \ \ \ \ \ \ \ \ \ \
\ \ \ \ \ \ \ \ \ \ \ \ \ \ \ \ \ \ \ \ \ \ \ \ \ \ \ \ \ \ \ \ \ \ \ \ \ \
\ \ \ \ \ \ \ \ \ \ Q.E.D.

\bigskip \ \ \ \ \ \ \

\ \ \ \ \ \ \ \ \ \ \ \ \ \ \ \ \ \ \ \ \ \ \ \ \ \ \ \ \ \ \ \ References \
\ \ \ \ \ \ \ \ \ \ \ \ \ \ \ \ \ \ \ \ \

\bigskip

[BdM] Boutet de Monvel, L., Int\'{e}gration des \'{e}quations de
Cauchy-Riemann induites formelles, S\'{e}minaire Goulaouic-Lions-Schwartz,
Expos\'{e} IX, 1974-1975.

\bigskip

[BG] Beals, R. and Greiner, P., Calculus on Heisenberg manifolds, Ann.
Math.
Studies, 119 (1988), Princeton University Press.

\bigskip

[Bo] Bogomolov, Fillability of contact pseudoconvex manifolds, G\"{o}ttingen
Univ. preprint. Heft 13 (1993), 1-13.

\bigskip
 
[Ca] Cartan, E., Sur la g\'{e}ometrie pseudo-conforme des hypersurfaces de
l'espace de deux variables complexe, I, Ann. Mat. 11 (1932), 17-90; II, Ann.
Sc. Norm. Sup. Pisa 1 (1932), 333-354.

\bigskip

[Ch] Cheng, J.-H., Contact topology and $CR$ geometry in three dimensions,
to appear in Proceedings of the 3rd Asian Mathematical Conference 2000 
(Quezon City, Philippines); math.SG/0106071.

\bigskip

[CL1] Cheng, J.-H. and Lee, J. M., The Burns-Epstein invariant and
deformation
of $CR$ structures, Duke Math. J. 60 (1990), 221-254.\ \ \ \ \ \ \ \ \ \ \ \
\ \

\bigskip

[CL2] ---------------, A local slice theorem for 3-dimensional $CR$
structures, Amer. J. Math. 117 (1995), 1249-1298.\ \ \ \ \ \ \ \ \ \ \ \ \ \
\ \ \ \ \ \ \ \ \ \ \ \ \ \

\bigskip

[El] Eliashberg, Y., Symplectic geometry of plurisubharmonic functions,
Notes by M. Abreu, NATO Adv. Sci. Inst. Ser. C, Math. Phys. Sci., 488,
Gauge theory and symplectic geometry (Montreal, PQ, 1995), 49-67,
Kluwer Acad. Publ., Dordrecht, 1997.

\bigskip
 
[Ko] Kohn, J. J., The range of the tangential Cauchy-Riemann operator, Duke
Math. J., 53 (1986), 525-545.

\bigskip

[L1] Lee, J. M., The Fefferman metric and pseudohermitian invariants, Trans.
Amer. Math. Soc., 296 (1986), 411-429. \ \ \ \ \ \ \ \ \ \ \ \ \ \ \ \ \ \ \
\ \ \ \ \ \ \ \ \ \ \ \ \ \ \ \ \ \ \ \ \ \ \ \ \ \ \ \ \ \ \ \

\bigskip

[L2] -----------, Pseudo-Einstein structures on $CR$ manifolds, Amer. J.
Math., 110 (1988), 157-178. \ \ \ \ \ \ \ \ \ \ \ \ \ \ \ \ \ \ \ \ \ \ \ \
\ \ \ \ \ \ \ \ \ \ \ \ \ \ \ \ \ \ \ \ \ \ \ \ \ \ \ \ \ \ \ \ \ \ \ \ \ \
\ \

\bigskip

[Lem] Lempert, L., On three-dimensional Cauchy-Riemann manifolds, J. Amer.
Math. Soc. 5 (1992), 923-969.

\bigskip

[NN] Newlander, A. and Nirenberg, L., Complex analytic coordinates in almost
complex manifolds, Ann. Math. 65 (1957), 391-404. \ \ \ \ \ \ \ \ \ \ \ \ \
\ \ \ \ \ \ \

\bigskip

[See] Seeley, R. T., Complex powers of an elliptic operator, Amer. Math.
Soc.
Proc. Symp. Pure Math. 10 (1967), 288-307.\

\bigskip

[Ta] Tanaka, N., A differential geometric study on strongly pseudo-convex
manifolds, 1975, Kinokuniya Co. Ltd., Tokyo. \ \ \ \ \ \ \ \ \ \ \ \ \ \ \ \
\ \

\bigskip

[Tr] Treves, F., Basic linear partial differential equations, Academic
Press,
New York (1975).

\bigskip

[We] Webster, S. M., Pseudohermitian structures on a real hypersurface, J.
Diff. Geom., 13 (1978), 25-41.

\bigskip

\bigskip

Jih-Hsin Cheng

Institute of Mathematics, Academia Sinica

Taipei, R.O.C. (Taiwan)

E-mail address: cheng@math.sinica.edu.tw

\end{document}